\documentclass[a4paper]{amsart}

\usepackage{amssymb,amsmath}
\usepackage{latexsym}
\usepackage{eucal}
\makeatletter
\@addtoreset{equation}{section}\makeatother

\newtheorem{theo}{Theorem}[section]
\newtheorem{lem}[theo]{Lemma}
\newtheorem{prop}[theo]{Proposition}

\newcommand{\mc}{\mathcal}
\newcommand{\rr}{\mathbb{R}}
\newcommand{\nn}{\mathbb{N}}
\newcommand{\cc}{\mathbb{C}}
\newcommand{\hh}{\mathbb{H}}

\newcommand{\la}{\lambda}
\newcommand{\eps}{\epsilon}

\newcommand{\pl}{\partial}
\newcommand{\x}{\times}

\newcommand{\til}{\widetilde}
\newcommand{\bbar}{\overline}

\newcommand{\supp}{\textrm{Supp}}
\newcommand{\cjd}{\rangle}
\newcommand{\cjg}{\langle}

\newcommand{\demi}{\frac{1}{2}}
\newcommand{\ndemi}{\frac{n}{2}}

\newcommand{\beq}{\begin{equation}}
\newcommand{\eeq}{\end{equation}}
\newcommand{\argth}{\textrm{argtanh}}

\newcommand{\indic}{\operatorname{1\negthinspace l}}
\def\qed{\hfill$\square$}
\begin{document}

\title[Absence of resonances near the critical line]
{Absence of resonance near the critical line
on asymptotically hyperbolic spaces}
\author[Colin Guillarmou]{Colin Guillarmou}
\address{Laboratoire de Math\'ematiques Jean Leray\\
         UMR 6629 CNRS/Universit\'e de Nantes \\
         2, rue de la Houssini\`ere \\
         BP 92208 \\
         44322 Nantes Cedex 03\\
         France}
     \email{cguillar@math.univ-nantes.fr}
\subjclass[2000]{Primary 58J50, Secondary 35P25}
\maketitle
\begin{abstract}
\noindent As a consequence of a result of Cardoso-Vodev, we show
that the resolvent of the Laplacian on asymptotically hyperbolic
manifolds is analytic in an exponential neighbourhood of the
critical line $\{\Re(\la)=\ndemi\}$. The case of non-trapping
metrics with constant curvature near infinity is also considered:
there is a strip $\{\Re(\la)>\ndemi-\eps\}$ with a finite number
of resonances.

\end{abstract}
\vspace{0.5cm}

\section{Introduction}

The purpose of this note is to give some `free of resonance' regions near the
critical line for the Laplacian on asymptotically hyperbolic manifolds.

An asymptotically hyperbolic manifold is a smooth non-compact
Riemannian manifold $(X,g)$ of dimension $n+1$ which is the
interior of a smooth compact manifold with boundary
$\bar{X}=X\cup\pl\bar{X}$ and such that for all boundary defining
function $x$ of $\bar{X}$ (i.e. $\pl\bar{X}=\{x=0\}$ and
$dx|_{\pl\bar{X}}\not=0$), $x^2g$ extends to a smooth metric on
$\bar{X}$ and $|dx|_{x^2g}=1$ on $\pl\bar{X}$. The metric can then
be expressed in a collar neighbourhood of the boundary
$(0,\eps)_x\x\pl\bar{X}_y$ by
\begin{equation}\label{modelform}
g=\frac{dx^2+h(x,y,dy)}{x^2}
\end{equation}
with $h(x,y,dy)$ a smooth tensor up to the boundary $\{x=0\}$.
It can be seen that $(X,g)$ is a complete manifold with
curvatures approaching $-1$ near the boundary $\pl\bar{X}$
(the boundary is the infinity of $\bar{X}$ with respect to the metric
$g$) and that the hyperbolic convex co-compact quotients
are contained in this class of manifolds.

It is well known that the spectrum of the Laplacian $\Delta_g$ acting on functions splits into
absolutely continuous spectrum $[\frac{n^2}{4},\infty)$ and a
finite set of eigenvalues $\sigma_{pp}(\Delta_g)\subset
(0,\frac{n^2}{4})$. The modified resolvent
\[R(\la):=(\Delta_g-\la(n-\la))^{-1}\]
is then meromorphic on $\{\Re(\la)>\ndemi\}$ with finite rank poles at
each $\la_e$ satisfying $\la_e(n-\la_e)\in\sigma_{pp}(\Delta_g)$.
Mazzeo and Melrose \cite{MM} have constructed the meromorphic extension of
$R(\la)$ to $\cc\setminus \demi(n-\nn)$ with poles of finite multiplicity,
which are called \textsl{resonances}. Physically, the most
interresting resonances are those in a neighbourhood of
the critical line $\{\Re(\la)=\ndemi\}$
(which corresponds to the essential spectrum).
It turns out that their localization is closely related
to the number of geodesics trapped in compact sets of $(X,g)$.
However, a general principle which seems to hold for geometric 
scattering on a large class of infinite volume manifolds
is that there exists a free of resonance 
region of the form
\[\{\la\in\cc; \Im(\la)\leq e^{-C_1|\la|}, |\Re(\la)|\geq C_2\},\quad C_1,C_2>0\]
when the critical line is the axis $\{\Im(\la)=0\}$.
This was first proved by Burq \cite{BU} in Euclidean scattering 
and Vodev \cite{V1} on some surfaces with negative constant curvature near
infinity.
It is worth noting that these results are optimal in the sense that, 
in general, the existence of elliptic closed geodesics implies the existence
of resonances which are exponentially close to the critical
line (see \cite{S,P1,P2}).

Here, we deal with the case of asymptotically hyperbolic manifolds:

\begin{theo}\label{main}
Let $(X,g)$ be an asymptotically hyperbolic manifold and $x$ a
boundary defining function,
then there exist $C_1,C_2>0$ such that
the weighted resolvent $x^\demi(\Delta_g-\la(n-\la))^{-1}x^\demi$
extends analytically from $\{\la\in\cc;|\Im(\la)|>C_2, \Re(\la)>\ndemi\}$ to
\begin{equation}\label{exponvois}
\{\la\in\cc; |\Im(\la)|>C_2,
\Re(\la)> \ndemi-e^{-C_1|\la|}\}.
\end{equation}
as bounded operators on $L^2(X):=L^2(X,dvol_g)$.
\end{theo}

The essential ingredients are a uniform bound of the weighted resolvent
norm on the critical line and a sharp parametrix of the
meromorphically continued resolvent.
In our case, the resolvent bound $||\rho R(\la)\rho||\leq Ce^{C|\la|}$
on the critical line has been proved by Cardoso and Vodev
\cite{CV}, $\rho$ being a
weight function decreasing to $0$ near infinity.
To extend analytically $\rho R(\la)\rho$ to the region (\ref{exponvois}),
the main point is to see it as a perturbation
of the resolvent of the Laplacian on a model space $X_0$ which
is sufficiently close to our manifold. A good candidate for $X_0$
is the warped product $(0,\eps)_x\x\pl\bar{X}_y$ equipped with the metric 
\[g_0:=\frac{dx^2+h(0,y,dy)}{x^2}\]
and take Dirichlet condition at $x=\eps$, but for technical
reasons we will better use $(0,\infty)\x \pl\bar{X}$ with the
same metric and localize the resolvent near $x=0$ with
cut-off functions (a similar approach is used by Vodev \cite{V1}).
This model resolvent $R_0(\la)$ needs to have an analytic
extension on weighted spaces
in a neighbourhood of the form (\ref{exponvois}), with a norm
bounded by $Ce^{C|\la|}$. The classical resolvent equation
\[R(\la)-R(z)=(\la(n-\la)-z(n-z))R(\la)R(z)\]
in the physical sheet $\{\Re(\la)>\ndemi\}$ and the
approximation of $R(\la)$ by $R_0(\la)$ allow to write
\[\rho R(\la)\rho (1+(\la-z)K(\la,z))=K_1(\la,z)\]
where $K(\la,z), K_1(\la,z)$ are some operators which are expressed
in terms of the
model resolvents, $\rho R(z)\rho$ and some error terms. At last,
the extension properties of $\rho R_0(\la)\rho$ through the
critical line and those of $\rho R(z)\rho$ up to this critical line can
be used to extend $K(\la,z)$ and $K_1(\la,z)$ to $z\in
\{\Re(\la)=\ndemi\}$ and $\la$ in the neighbourhood
(\ref{exponvois}). The bound on the norm of $\rho R(z)\rho$ and
$\rho R_0(\la)\rho$ can then be used to show that
$|\la-z|.||K(\la,z)||\leq \demi$ for $\Re(z)=\ndemi$ and
$|\la-z|\leq C^{-1}e^{-C|z|}$ with $C>0$ large and
independent of $(\la,z)$; this allows to invert holomorphically
$1+(\la-z)K(\la,z)$ and to define $\rho R(\la)\rho$ in
(\ref{exponvois}).

In the case of non-trapping metrics, Vodev \cite{V2}
proved that the norm of the resolvent on the critical
line grows not faster than $C|\la|^{-1}$ when $|\Im(\la)|\to \infty$,
which implies a larger extension of the resolvent through the essential
spectrum. We especially consider the case of non-trapping manifolds
with constant curvature near infinity and obtain:

\begin{theo}\label{hyper}
Let $(X,g)$ be a conformally compact manifold with constant curvature
outside a compact subset and let $x$ be a boundary defining function. If $g$ is
non trapping, there exist $C_1,C_2>0$ such that the weighted
resolvent $x^\demi(\Delta_g-\la(n-\la))^{-1}x^\demi$ extends
analytically from $\{\la\in\cc;|\Im(\la)|>C_2, \Re(\la)>\ndemi\}$ to
\[\{\la\in\cc; |\Im(\la)|>C_2, \Re(\la)> \ndemi-C_1\}\]
as bounded operators on $L^2(X)$.
\end{theo}

Note that this non-trapping condition is not
satisfied for non-elementary convex co-compact
quotients of $\hh^{n+1}$. When $X=\Gamma\backslash\hh^{n+1}$
with $\Gamma$ a non-elementary convex co-compact group of isometries,
better results are available with the help of Selberg's zeta function:
it is well known that there exists a half plane $\Re(\la)>\delta$
with no resonance, where $\delta$ is the dimension of the limit set
(see also a result of Naud \cite{N} in dimension $2$).
This shows that Theorem \ref{hyper} is weak in the sense
that we need the non-trapping assumption but the compensation is that
we do not have to confine ourselves to the rigid class of constant
curvature manifolds.\\

The paper is organized as follows: in Section 1, we recall
Cardoso-Vodev Theorem, then we study our models in Section 2
and finally we give the proof of the results in Section 3.\\

\textbf{Acknowledgements.} I would like to thank Georgi Vodev who
pointed out to me this problem.

\section{Cardoso-Vodev result}

In \cite{CV}, Cardoso and Vodev consider some Riemannian manifolds $X$
with controlled structure near infinity
and they obtain exponential bounds for the weighted resolvent norm on the critical
line. These manifolds have the following properties outside a compact set $Z$
\begin{equation}\label{decomposition}
X\setminus Z \cong \left( [R,\infty)\x
S, g:=dr^2+\sigma(r)\right), \quad R\gg 1,
\end{equation}
where $S$ is a n-dimensional smooth compact manifold, $\cong$
means `isometric' and
$\sigma(r)=\sigma(r,y,dy)$ is a family of metrics on
$S_r:=\{r\}\x S$ which satisfy
\begin{equation}\label{propsigma}
|q(r,y)|\leq C ,\quad
\pl_r q(r,y)\leq Cr^{-1-\delta}, \quad \delta>0,\quad r>R,
\end{equation}
\begin{equation}\label{prop2sigma}
-\pl_r(\sigma^{-1})(r,y,\xi)\geq\frac{C}{r}\sigma^{-1}(r,y,\xi)
,\quad \forall (y,\xi)\in T^*S_r,
\end{equation}
$\sigma^{-1}(r)$ being the
principal symbol of the Laplacian on $(S_r,\sigma(r))$
and $q(r,y)$ is an effective potential defined by
\[q(r,y):=(2^{-1}\pl_r\log\nu)^2+(2\nu)^{-2}\sum_{i,j}\sigma^{ij}\pl_{y_i}\nu
\pl_{y_j}\nu+2^{-1}\nu\Delta_g(\nu^{-1})\]
with $\nu:=(\det(\sigma_{ij}))^{\demi}$. An asymptotically
hyperbolic metric $g=x^{-2}(dx^2+h(x,y,dy))$
can be decomposed as in (\ref{decomposition}) by
putting $x=e^{-r}$ and we get
\[\sigma(r,y,dy)=e^{2r}h(e^{-r},y,dy)\]
with $h(x,y,dy)$ smooth up to $x=0$. In $x$ coordinate, we have
$\nu\in x^{-n}C^\infty(\bar{X})$, $\pl_r=-x\pl_x$ and
$h^{ij}\in C^\infty(\bar{X})$ thus
\[(x\pl_x\log\nu)^2+\nu^{-2}\sum_{i,j}x^2h^{ij}\pl_{y_i}\nu
\pl_{y_j}\nu \in C^\infty(\bar{X})\]
\[\nu\Delta_g\nu^{-1}=\nu\Big(-(x\pl_x)^2-(x\pl_x\log\nu)x\pl_x-
\nu^{-1}\sum_{i,j}\pl_{y_i}(\nu x^2h^{ij}\pl_{y_j})\Big)\nu^{-1}\in
C^\infty(\bar{X})\]
and (\ref{propsigma}) is satisfied for all $\delta>0$. Moreover
we have for all $(y,\xi)\in T^*\pl\bar{X}$
\[\frac{x\pl_x (x^2h^{-1})(x,y,\xi)}{x^2h^{-1}(x,y,\xi)}=
2+x\frac{\pl_x(h^{-1})(x,y,\xi)}{h^{-1}(x,y,\xi)}\geq 1\]
if $x\leq \eps$ with $\eps$ small, and we obtain that
(\ref{prop2sigma}) is satisfied. As a conclusion, asymptotically
hyperbolic manifolds are in the class of manifolds studied by
Cardoso and Vodev \cite{CV} and Vodev \cite{V2}, so their results can
be summarized in that case in the
\begin{theo}\label{cardvod}
Let $(X,g)$ be an asymptotically hyperbolic manifold and $x$
a boundary defining function. There exists
$C>0$ such that the weighted resolvent
$x^\demi(\Delta_g-\la(n-\la))^{-1}x^\demi$
extends continuously from $\{\Re(\la)>\ndemi, |\Im(\la)|\geq 1\}$
to $\{\Re(\la)\geq \ndemi, |\Im(\la)|\geq 1\}$ on $L^2(X)$
and the extension satisfies 
\[||x^{\demi}R(\la)x^{\demi}||_{\mc{L}(L^2,H^p)}\leq Ce^{C|\la|},\quad C>0\]
for $p=0,1$, $|\Im(\la)|\geq 1$ and $0\leq\Re(\la)-\ndemi\leq 1$, 
where $H^p$ means the p-Sobolev space on $X$ with respect to the metric $g$.
If in addition $g$ is non-trapping we have for $p=0,1$,
$|\Im(\la)|\geq 1$ and $0\leq\Re(\la)-\ndemi\leq 1$
\[||x^{\demi}R(\la)x^{\demi}||_{\mc{L}(L^2,H^p)}\leq C|\Im(\la)|^{-1+p},\quad
C>0.\]
\end{theo}

\section{Two models}

Before giving the models, we recall a few properties of some
differential operators on $X$. Let $(X,g)$ be a conformally compact
manifold and $\Delta_g$ the Riemannian Laplacian. If $x$ is a boundary
defining function and $(y_i)_{i=1,\dots,n}$ some coordinates on
$\pl\bar{X}$, the space $\mc{V}_0(\bar{X})$ of smooth vector
fields on $\bar{X}$ which vanish on $\pl\bar{X}$ is locally
generated by $x\pl_x,x\pl_{y_i}$ for $i=1,\dots,n$ near the
boundary. We denote by $\textrm{Diff}_0^k(\bar{X})$ the space of differential
operators of order $k$ generated by $k$ products of elements of
$\mc{V}_0(\bar{X})$
\[\textrm{Diff}_0^k(\bar{X}):=\underset{0\leq i\leq k}{\textrm{Vect}}
\mc{V}_0(\bar{X})^k,\quad  \mc{V}_0(\bar{X})^0=C^\infty(\bar{X}).\]
For example it is straightforward to check that
\[\Delta_g\in \textrm{Diff}_0^2(\bar{X}).\]
$\Delta_g$ is now considered as the self-adjoint operator obtained
by Friedrichs extension from the Laplacian on $C_0^\infty(X)\subset
L^2(X)=L^2(X,dvol_g)$. For $k\in\rr$, we define the k-Sobolev space by
\[H^k(X):=\mc{D}om((1+\Delta_g)^{\frac{k}{2}}),\]
where $\mc{D}om$ means the domain.
The Sobolev spaces associated to two different
conformally compact metrics are the same
(for instance, it is done for $k=1,2$ by Froese-Hislop
\cite[appendix]{FH} in a more general framework) and
\begin{equation}\label{propdiff0}
\forall D^k\in \textrm{Diff}^k_0(\bar{X}), \quad D^k\in
\mc{L}(H^s(X),H^{s-k}(X)),\quad s\in [0,2].
\end{equation}
Moreover a useful
property of these differential operators is the following
\begin{equation}\label{commutdiff0}
x^{-\alpha}D^kx^\alpha\in \textrm{Diff}_0^k(\bar{X}), \quad
\alpha\in\rr, \quad D^k\in\textrm{Diff}^k_0(\bar{X})
\end{equation}
which is easily seen from the commutator $[x\pl_x,x^\alpha]$
in local charts near $\pl\bar{X}$.\\

Let us now study two models which will be respectively
used for the parametrix construction of the general
case and for the case of constant curvature near infinity.\\

Let $(M,h_0)$ be a Riemannian compact manifold of dimension $n$ and
\begin{equation}\label{model}
X_0:=(0,+\infty)_x\x M, \quad g_0:=x^{-2}(dx^2+h_0).
\end{equation}
Though $(X_0,g_0)$ is not conformally
compact (there is a cusp end when $x\to\infty$),
it has a conformally compact structure near $x=0$. We
could take as model operator the Laplacian on
$((0,1]\x\pl\bar{X},g_0)$ with Dirichlet condition at $x=1$ to
have a conformally compact structure (with boundary), but
we prefer to use $X_0$ since it carries more symmetry and is
therefore easier to study.

As for conformally compact manifolds,
let us denote $\bar{X}_0:=[0,+\infty) \x M$ and
$\textrm{Diff}^k_0(\bar{X}_0)$
the space of smooth differential operators of order $k$
on $\bar{X}_0$ with support in $[0,1]\x M$ and which can be locally
written
\[\sum_{i+|\alpha|\leq k}a_{i,\alpha}(x,y)(x\pl_x)^ix^{|\alpha|}
\pl^\alpha_{y_i}, \quad a_{i,\alpha}\in C^\infty(\bar{X}_0),\]
where $(y_i)_{i=1,\dots,n}$ are some local coordinates on $M$.

By taking the new variable $r=\log x$, it is easy to see that the
Laplacian $\Delta_{g_0}$ is unitarily equivalent to
\[P_0=-\pl_r^2+e^{2r}\Delta_{h_0}+\frac{n^2}{4}\]
on $L^2(\rr\x M, drdvol_{h_0})$.
As before, we define the Sobolev spaces by
\[H^k(X_0):=\mc{D}om((1+\Delta_{g_0})^{\frac{k}{2}})\cong
\mc{D}om((1+P_0)^{\frac{k}{2}}).\]
We first remark that
the arguments given by Froese and Hislop
in \cite[appendix]{FH} prove that for $k=0,1,2$ and $s\in[0,2]$
\begin{equation}\label{propdiff0x}
\forall D^k\in \textrm{Diff}^k_0(\bar{X}_0),\quad
D^k\in \mc{L}(H^{s}(X_0),H^{s-k}(X_0))
\end{equation}
though they do not consider the `cuspidal' part $\{r\in \rr^+\}$.
Of course the case $k=1$ can be directly
obtained from the identity
\begin{equation}\label{sobolevh1}
||(1+P_0)^\demi u||^2_{L^2}=
||\pl_ru||_{L^2}^2+||e^r\Delta^\demi_{h_0}u||^2_{L^2}
+\left(\frac{n^2}{4}+1\right)||u||^2_{L^2}.
\end{equation}

We shall first see how the resolvent of $\Delta_{g_0}$ can be
extended to the non-physical sheet and we will give an upper bound of
its weighted norm.

\begin{lem}\label{casmodel}
Let $x_0<1$, $(X_0,g_0)$ defined in (\ref{model})
and $\rho=\rho(x)$ a smooth function on $X_0$ with support in $\{x<1\}$
such that $\rho(x)=x^\demi$ for $x\leq x_0$. Then the weighted
resolvent
\[\rho R_0(\la)\rho:=\rho(\Delta_{g_0}-\la(n-\la))^{-1}\rho\]
extends analytically from $\{\Re(\la)>\ndemi,|\Im(\la)|>1\}$
to $\{\Re(\la)>\ndemi-\frac{1}{4}, |\Im(\la)|> 1\}$ and it satisfies
\begin{equation}\label{h1bound}
||\pl_\la^q \rho R_0(\la)\rho||_{\mc{L}(L^2,H^p)}\leq
e^{C|\la|}, \quad C>0
\end{equation}
in that region for $q=0,1$, $p=0,1,2$.
\end{lem}
\textsl{Proof}: let us begin by $p=2$. Let $\chi\in C_0^\infty([0,1))$
which is equal to $1$ on $\supp(\rho)$ and $\chi_1\in C_0^\infty([0,1))$
such that $\chi_1=1$ on $\supp(\chi)$. We then have for
$\la_0\in \{\Re(\la)>\ndemi\}$
\[
(\Delta_{g_0}+1)\rho R_0(\la)\rho=
[\Delta_{g_0},\rho]x^{-\demi}\chi(x)x^\demi R_0(\la)\rho+\rho^2+(\Lambda+1)\rho
R_0(\la)\rho,
\]
\[
\chi(x)x^\demi R_0(\la)\rho=R_0(\la_0)\Big([\Delta_{g_0},\chi(x) x^\demi]x^{-\demi}
\chi_1(x)x^\demi R_0(\la)+x^\demi+(\Lambda-\Lambda_0)\chi(x) x^\demi
R_0(\la)\Big)\rho
\]
where $\Lambda:=\la(n-\la),\Lambda_0:=\la_0(n-\la_0)$.
Since $[\Delta_{g_0},\rho]x^{-\demi}\in\textrm{Diff}^1_0(\bar{X}_0)$ and
$[\Delta_{g_0},\chi(x)x^\demi]x^{-\demi}\in\textrm{Diff}^1_0(\bar{X}_0)$
in view of (\ref{commutdiff0}), it suffices to prove (\ref{h1bound})
for $p=0$ and use (\ref{propdiff0x}) to obtain the other cases.

Let us now use $P_0$ instead of $\Delta_{g_0}$. 
We have a decomposition induced by the spectral resolution of $\Delta_{h_0}$
\[P_0=\bigoplus_{j\in\nn_0} P^{(j)}_0, \quad
P^{(j)}_0:=-\pl^2_r+e^{2r}\mu^2_j+\frac{n^2}{4}\]
where $(\mu_j^2)_{j\in\nn_0}$ are the eigenvalues of $\Delta_{h_0}$
(counted with multiplicities) associated to
an orthonormal basis of $L^2(M)$ of eigenvectors $(\psi_j)_{j\in\nn_0}$.
If we denote the resolvent of $P^{(j)}_{0}$ on $L^2(\rr,dr)$ by
\[R^{(j)}_{0}(\la):=(P^{(j)}_{0}-\la(n-\la))^{-1}\]
we clearly have for $f\in L^2(X)$ and $f_j:=\cjg
f,\psi_j\cjd_{L^2(M)}$
\begin{equation}\label{formuler}
\rho(\Delta_{g_0}-\la(n-\la))^{-1}\rho f= \sum_{j\in\nn_0}(\rho
R^{(j)}_{0}(\la)\rho f_j)\psi_j.
\end{equation}
Note that for $\mu_j\not=0$ the translation
\[U_j:\left\{
\begin{array}{rcl}
L^2(\rr,dr)&\to& L^2(\rr,dr)\\
f&\to& f(\log(\mu_j)+\bullet)
\end{array}\right.\]
is an isometry and that $U_j^{-1}P^{(j)}_0U_j=Q$ with
$Q:=-\pl^2_r+e^{2r}+\frac{n^2}{4}$. Let us set $k=\la-\ndemi$ for
simplicity. The Green kernel for $Q$ is then easy to find for
$\Re(k)>0$ (see \cite[Ex. 4.15]{T})
\[R_Q(\la;r,t):=(Q-\la(n-\la))^{-1}(r,t)=
-K_{-k}(e^r)I_k(e^t)H(r-t)-I_k(e^r)K_{-k}(e^t)H(t-r)\]
with $H$ the Heaviside function and $I_k,K_k$ the modified Bessel functions
whose integral representations (when they converge) are
\begin{equation}\label{besseli1}
I_k(z)=\frac{2^{1-k}z^k}{\Gamma(k+\demi)\Gamma(\demi)}
\int_0^1(1-u^2)^{k-\demi}\cosh(zu)du
\end{equation}
\begin{equation}\label{besseli2}
I_k(z)=\frac{1}{\pi}
\int_0^\pi e^{z\cos(u)}\cos(ku)du-\frac{\sin(k\pi)}{\pi}\int_0^\infty 
e^{-z\cosh(u)-ku}du
\end{equation}
\begin{equation}\label{besselk1}
K_{k}(z)=K_{-k}(z)=\frac{\Gamma(-k+\frac{3}{2})2^{-k+1}}{\Gamma(\demi)}z^{k}
\int_0^\infty \frac{\sin(tz)}{tz}t^2(1+t^2)^{k-\frac{3}{2}}dt\end{equation}
\begin{equation}\label{besselk2}
K_{-k}(z)=\int_0^\infty\cosh(ku)e^{-z\cosh(u)}du.
\end{equation}
Moreover for $\mu_j=0$, the expression of the meromorphically extended
euclidian resolvent kernel is well-known in $\cc\setminus\{0\}$
\[R^{(j)}_0(\la;r,t)=(2k)^{-1}e^{-k|r-t|}\]
and we have for $\Re(\la)>\ndemi-\frac{1}{4}$, $|\Im(\la)|\geq 1$ and
$p=0,1,2$ 
\begin{equation}\label{mujegal0} 
||\pl_r^p \rho(e^\bullet)
R^{(j)}_0(\la)\rho(e^\bullet)||_{\mc{L}(L^2(\rr))}\leq
C\left|\la-\ndemi\right|^{-1+p} \textrm{ if } \mu_j=0. \end{equation} 
From
these expressions, one can remark that there is no resonance
except $\ndemi$ for this problem. Without loss of generality, we can
suppose that $\rho(e^r)=e^{\frac{r}{2}}\chi(r)$ where $\chi$ is a smooth
function on $\rr$ such that $\chi(r)=1$ when $r\leq -1$ and
$\chi(r)=0$ when $r\geq 0$.
Since from (\ref{besseli1})-(\ref{besselk2}) we have
\[|I_k(e^t)\rho(e^{t-\log(\mu_j)})|\leq  
\left\{\begin{array}{l}
e^{C|k|}e^{e^t}\textrm{ for }t>0\\
 e^{C|k|}e^{\frac{t}{4}}\textrm{ for }t\leq 0
\end{array}\right.\]
and
\[|K_{-k}(e^r)\rho(e^{r-\log(\mu_j)})|\leq 
\left\{\begin{array}{l}
e^{C|k|}e^{-e^r} \textrm{ for }r>0\\
e^{C|k|}e^{\frac{r}{4}} \textrm{ for }r\leq 0
\end{array}\right.\]
for $|\Re(k)|\leq \frac{1}{4}$ and $|\Im(k)|\geq 1$, we can easily deduce that
\[|K_{-k}(e^r)\rho(e^{r-\log(\mu_j)})|
\int_{-\infty}^r|I_k(e^t)\rho(e^{t-\log(\mu_j)})|dt\leq e^{C|k|}\]
\[|I_{k}(e^r)\rho(e^{r-\log(\mu_j)})|
\int_{r}^\infty|K_{-k}(e^t)\rho(e^{t-\log(\mu_j)})|dt\leq e^{C|k|}\]
and Schur's lemma implies that
\[||U_j^{-1}(\rho)R_Q(\la)U_j^{-1}(\rho)||_{\mc{L}(L^2)}\leq e^{C|k|}\]
for $|\Re(k)|\leq \frac{1}{4}$, $|\Im(k)|\geq 1$. But since
$U_j^{-1}(\rho(e^{\bullet}))=U_j^{-1}\rho(e^\bullet) U_j$ as 
operators and $U_j$ is an isometry, we use
\[R_0^{(j)}(\la)=U_jR_Q(\la)U_j^{-1}\]
to conclude that for $|\Im(\la)|\geq 1$ and
$|\Re(\la)-\ndemi|\leq \frac{1}{4}$
\begin{equation}\label{borneq} 
||\rho(e^\bullet) R^{(j)}_0(\la)\rho(e^{\bullet})||_{\mc{L}(L^2)}\leq e^{C|\la|}.
\end{equation}
Finally we combine (\ref{borneq}) with (\ref{mujegal0}) and (\ref{formuler}).
The bound for $q=1$ (one derivative with respect to $\la$) is directly
obtained from the case $q=0$ and Cauchy's formula.
\qed\\

\textsl{Remark}: a better estimate can be obtained
but we do not need it for our purpose.\\

The second model is the hyperbolic space $(\hh^{n+1},g_h)$ with
its usual metric. To see $\hh^{n+1}$ as a conformally compact
manifold we take the Beltrami model
\[\hh^{n+1}=\{m\in\rr^{n+1}; |m|<1\}, \quad g_h=\frac{4dm^2}{(1-|m|^2)^2}\]
and we set for example $x=2(1-|m|)(1+|m|)^{-1}$
as boundary defining function.

\begin{prop}\label{bound:uniform}
The weighted hyperbolic resolvent
\[x^\demi R_h(\la)x^\demi:=x^\demi(\Delta_{g_h}-\la(n-\la))^{-1}x^\demi\]
extends holomophically from $\{\Re(\la)>\ndemi\}$ to
$\{\Re(\la)>\ndemi-\frac{1}{8}\}$ with values 
in $\mc{L}(L^2(\hh^{n+1}),H^1(\hh^{n+1}))$ and it satisfies
\begin{equation}\label{normresolvlib}
\left|\left|\pl_\la^q x^\demi R_h(\la)x^\demi\right|\right|_
{\mc{L}(L^2,H^p)}\leq C\left|\la-\ndemi\right|^{-1+p},\quad C>0
\end{equation}
in the same region (with $\la\not=\ndemi$) for $p=0,1$, $q=0,1$.
\end{prop}
\textsl{Proof}: to show these bounds, we consider the wave
operators
\[U_0(t):=\cos\left (t\sqrt{\Delta_{g_h}-\frac{n^2}{4}}\right ),
\quad U_1(t):=\left (\Delta_{g_h}-\frac{n^2}{4}\right )^{-\demi}
\sin\left (t\sqrt{\Delta_{g_h}-\frac{n^2}{4}}\right ).\]
Let $\delta>0$ be small, $\chi\in C_0^{\infty}(\rr^+)$ with support in
$[0,1-\frac{\delta}{2}[$ and which is equal to $1$ on
$[0,1-\delta]$. Then set
\[\chi_t(m):=\chi\left (\frac{4\argth |m|}{t}\right )\]
whose support is included in the hyperbolic ball of radius
$\frac{t}{2}$ and which is equal to $1$ in the hyperbolic ball of radius
$\frac{t(1-\delta)}{2}$.

When $n+1$ is odd, we have by Huygens principle (see \cite{H})
\begin{equation}\label{huygens}
\chi_tU_0(t)\chi_t=0, \quad t>0.
\end{equation}
Now we split $x^\demi$ in the following way
\begin{equation}\label{decomp1}
x^\demi U_0(t)x^\demi=(1-\chi_t)x^\demi U_0(t)(1-\chi_t)x^\demi+
(1-\chi_t)x^\demi U_0(t)\chi_tx^\demi+\chi_t x^\demi U_0(t)(1-\chi_t)x^\demi
\end{equation}
for $t>0$. It is clear that
\begin{equation}\label{norm1}
||\chi_tx^\demi||\leq C, \quad
||U_0(t)||\leq 1
\end{equation}
and remark that
$m\in\supp(1-\chi_t)$ only if $|m|\geq
\tanh\left(\frac{(1-\delta)t}{4}\right )$, that is
\[x=2\frac{1-|m|}{1+|m|}\leq
2\frac{1-\tanh(\frac{t(1-\delta)}{4})}{1+\tanh(\frac{t(1-\delta)}{4})}
=2e^{-\frac{t(1-\delta)}{2}}.\]
So we find
\begin{equation}\label{norm2}
\left|\left|(1-\chi_t)x^\demi\right|\right|\leq Ce^{-\frac{t(1-\delta)}{4}}.
\end{equation}
From (\ref{decomp1}), (\ref{norm1}) and (\ref{norm2}) we deduce that
\[\left|\left|x^\demi U_0(t)x^\demi\right|\right|\leq Ce^{-\frac{t(1-\delta)}{4}}.\]
It remains to use the Laplace transform of $U_0$
\[ \left
(\la-\ndemi\right )x^\demi R_h(\la)x^\demi= \int_0^{\infty}
e^{t(\ndemi-\la)}x^\demi U_0(t)x^\demi dt\]
and (\ref{normresolvlib}) is
proved when $p=0$ and $n+1$ odd by taking $\delta<\demi$.
To deal with the case $n+1$ even, we study
$x^\demi\pl_t U_1(t)x^\demi$ and use $U_0(t)=\pl_t U_1(t)$.
We have
\begin{equation}\label{decomp2}
\chi_t \pl_tU_1(t)\chi_t=\pl_t\left (\chi_t U_1(t)\chi_t\right )
-(\pl_t\chi_t)U_1(t)\chi_t-\chi_t U_1(t)(\pl_t\chi_t)
\end{equation}
and the Schwartz kernel of $U_1(t)$ is (see \cite{H})
\begin{equation}\label{u1}
U_1(t;x,y)=
C_n\left (\sinh^2\Big (\frac{t}{2}\Big )- \sinh^2\Big
(\frac{d_{\hh^{n+1}}(x,y)}{2}\Big )\right )_+^{-\ndemi}
\end{equation}
with $C_n\in\rr$. Hence, by construction
of $\chi_t$, the operators in (\ref{decomp2}) have a smooth kernel
with support in $\{(x,y)\in\hh^{n+1}\x\hh^{n+1},d_{\hh^{n+1}}(x,y)
<t(1-\frac{\delta}{2})\}$ and (\ref{u1}) implies that
there exists $T>0$ such that
\begin{equation}\label{normonde2}
\left|\left|\chi_t \pl_t U_1(t)\chi_t\right|\right|\leq
Ce^{-\frac{nt}{2}}
\end{equation}
for $t\geq T$. At last, we proceed as in
the odd case: we split $x^\demi$ on the support of $\chi_t$ and
outside, which shows (\ref{normresolvlib}) for $p=0$.

Now for $p=1$, it suffices to show that 
$Dx^\demi R_h(\la)x^\demi$ extends to 
$\{\Re(\la)>\ndemi-\frac{1}{8}\}$ in $\mc{L}(L^2)$
for a finite number of $D\in\textrm{Diff}^1_0(\bbar{\hh}^{n+1})$
which span $\textrm{Diff}^1_0(\bbar{\hh}^{n+1})$ over 
$C^\infty(\bbar{\hh}^{n+1})$. 
In fact since $Dx^\demi\in x^\demi\textrm{Diff}_0^1(\bbar{\hh}^{n+1})$, 
we will study $x^\demi DR_0(\la)x^\demi$. 
With $D$ as before, we have for $\Re(\la)>\ndemi$
\begin{equation}\label{transflap2}
x^\demi DR_0(\la)x^\demi=\int_0^\infty
e^{t(\ndemi-\la)}x^\demi DU_1(t)x^\demi dt
\end{equation}
but $||DU_1(t)||\leq C$ in view of the continuity of $D$ 
from $H^1(\hh^{n+1})$ to $L^2(\hh^{n+1})$. Moreover
\[\chi_tDU_1(t)\chi_t=D\chi_tU_1(t)\chi_t-[D,\chi_t]U_1(t)\chi_t\]
has a smooth kernel with compact support and it is 
straightforward to check (from (\ref{u1})) that
there exists $T>0$ such that for all $t\geq T$
\begin{equation}\label{normonde3}
||\chi_tDU_1(t)\chi_t||\leq C e^{-\frac{nt}{2}} \end{equation}
Splitting $x^\demi DU_1(t)x^\demi$ in the same way as
(\ref{decomp1}) and using (\ref{normonde3}), (\ref{transflap2}), one deduces
the bound (\ref{normresolvlib}) for $p=1$.
The case $q=1$ is obtained by the Cauchy formula from $q=0$.
\qed\\

\section{Parametrix construction and proof of the main result}

In this section, we will give the construction of a
parametrix for the resolvent $R(\la)$ of $\Delta_g$ on  an arbitrary
asymptotically hyperbolic manifold $(X,g)$ whose metric is
\[g=x^{-2}(dx^2+h(x))\]
in a collar $(0,\delta)_x\x \pl\bar{X}$ near the boundary,
$h(x)$ being a smooth family of metrics on $\pl\bar{X}$.
Of course this parametrix will only be sufficient
to approach the resolvent in
$\{\Re(\la)>\ndemi-\delta\}$ with $\delta>0$ small.
Roughly, the singularities at $x=0$ and $x'=0$ of the kernel of $R(\la)$
are controlled by those of the resolvent kernel of the
Laplacian induced by the model metric $x^{-2}(dx^2+h(0))$
near the boundary, already studied in the previous section.\\

\textsl{Proof of Theorem \ref{main}}: consider $(X,g)$ an
asymptotically hyperbolic manifold and $R(\la):=(\Delta_g-\la(n-\la))^{-1}$
the resolvent for $\Delta_g$ in the physical sheet
$\{\Re(\la)>\ndemi\}$. Let $V$ be an open subset
in $(X,g)$ isometric to $U:=(0,\delta)_x\x\pl\bar{X}$
(with $\delta>0$) equipped
with the metric $x^{-2}(dx^2+h(x))$ and $i:V\to U$ this
isometry. Note that it is always possible to renormalize $x$ (take
$x'=x\delta^{-1}$) to have $\delta=1$ and it does not change the
structure of the metric, so we suppose $\delta=1$.
We now consider $(X_0,g_0)$ the Riemannian manifold defined in
(\ref{model}) with
\begin{equation}\label{defx0g0}
(M,h_0):=(\pl\bar{X},h(0)).
\end{equation}
Let $\mc{I}_U$ and $\mc{R}_U$ be the following bounded operators
\[\mc{R}_U :\left\{\begin{array}{rcl}
L^2(X_0,dvol_{g_0})&\to& L^2(U,dvol_{g_0})\\
f&\to& f(\iota_U(\bullet))
\end{array}\right.\]
\[\mc{I}_U :\left\{\begin{array}{rcl}
L^2(U,dvol_{g_0})&\to& L^2(X_0,dvol_{g_0})\\
f&\to& \indic_Uf
\end{array}\right.\]
with $\iota_U$ the inclusion $U\subset X_0$ and $\indic_U$ the characteristic
function of $U$. Similarly one can define the operators
$\mc{I}_V$ and $\mc{R}_V$ induced by the inclusion $V\subset X$.
Since $i^*g_0$ and $g$ are quasi-isometric on $V$, we obtain that
$i^*:L^2(U,dvol_{g_0})\to L^2(V,dvol_{g})$ and
$i_*: L^2(V,dvol_{g})\to L^2(U,dvol_{g_0})$ are bounded.
We then set
\[I^*:=\mc{I}_Vi^*\mc{R}_U \in \mc{L}(L^2(X_0,dvol_{g_0}),
L^2(X,dvol_g)),\]
\[I_*:=\mc{I}_Ui_*\mc{R}_V \in\mc{L}(L^2(X,dvol_g),
L^2(X_0,dvol_{g_0})).\]

For $j=1,2,3,4$, let $\psi_j$ a smooth function on $\rr^+$
which is equal to $1$ in $[0,\frac{j}{5}]$ and
to $0$ in $[\frac{j+1}{5},+\infty)$;
$\psi_j$ will also be considered as
a function on $U$ depending only on $x$. With $\til{\psi}_j:=i^*\psi_j$,
we have as operators
\[I_*I^*\psi_j=\psi_j,\quad I^*I_*\til{\psi}_j=\til{\psi}_j.\]
It is easy to check that there exist
$D_R,D_L\in\textrm{Diff}^2_0(\bar{X})$ such that
\[\Delta_g-I^*\Delta_{g_0}\psi_3I_*= xD_R,
\quad \Delta_g-I^*\psi_3\Delta_{g_0}I_*=xD_L.\]
Firstly, we take $\Re(\la)>\ndemi$ and
$R_0(\la)$ is the resolvent considered in Lemma \ref{casmodel}.
Observe that
\begin{equation}\label{restechi1}
(\Delta_{g_0}\psi_3-\la(n-\la))\psi_2R_0(\la)\psi_1=\psi_1+
[\Delta_{g_0},\psi_2]R_0(\la)\psi_1
\end{equation}
since $\psi_3\psi_2=\psi_2$ and $\psi_2\psi_1=\psi_1$.
Let $\chi_1:=1-\til{\psi}_1$
and let $\chi_0$ be a smooth function with compact support on $X$ which
is equal to $1$ on the support of $\chi_1$.
For $\la_0\in \{\Re(\la)>\ndemi\}$ fixed,
$\Lambda:=\la(n-\la)$ and $\Lambda_0:=\la_0(n-\la_0)$ we set
\[R_{0R}(\la):=I^*\psi_2R_0(\la)\psi_1 I_*,\quad
E_R(\la):=R_{0R}(\la)+\chi_0R(\la_0)\chi_1,\]
\[L_R(\la):=[\Delta_g,\chi_0]R(\la_0)\chi_1+
(\Lambda_0-\Lambda)\chi_0R(\la_0)\chi_1+I^*[\Delta_{g_0},\psi_2]R_0(\la)\psi_1I_*
+xD_RR_{0R}(\la),\]
and we deduce from (\ref{restechi1})
\begin{equation}\label{inverseright}
(\Delta_g-\Lambda)E_R(\la)=1+L_R(\la).
\end{equation}
Similarly, we give a left parametrix for $\Delta_g-\Lambda$
\[R_{0L}(\la):=I^*\psi_1R_0(\la)\psi_2I_*,\quad
E_L(\la)=R_{0L}(\la)+\chi_1R(\la_0)\chi_0,\]
\[L_L(\la)=R_{0L}(\la)xD_L+I^*\psi_1R_0(\la)[\psi_2,\Delta_{g_0}]I_*+
\chi_1R(\la_0)[\Delta_g,\chi_0]+
(\Lambda_0-\Lambda)\chi_1R(\la_0)\chi_0,\]
and we have
\begin{equation}\label{inverseleft}
E_L(\la)(\Delta_g-\Lambda)=1+L_L(\la).
\end{equation}
Let $z,\la\in\{s\in\cc;\Re(s)>\ndemi\}$ and $Z:=z(n-z)$, from
(\ref{inverseright}) and (\ref{inverseleft}) we then obtain
\begin{equation}\label{egalpoids}
R(\la)=E_R(\la)-R(\la)L_R(\la),\quad
R(z)=E_L(z)-L_L(z)R(z).
\end{equation}
On the other hand, the resolvent equation
\[R(\la)-R(z)=(\Lambda-Z)R(\la)R(z)\]
combined with the first identity of (\ref{egalpoids}) yield
\begin{equation}\label{paramfinal}
x^\demi R(\la)x^\demi(1+K(\la,z))=K_1(\la,z),
\end{equation}
\[K(\la,z):=(\Lambda-Z)x^{-\demi}L_R(\la)R(z)x^\demi,\quad
K_1(\la,z):=x^\demi R(z)x^\demi+(\Lambda-Z)x^\demi
E_R(\la)R(z)x^\demi.\]
For simplicity, we will now denote by $D^p$ (resp. $D_0^p$)
all differential operator in $\textrm{Diff}^p_0(\bar{X})$
with support in $\supp(\til{\psi}_3)$
(resp. $\textrm{Diff}_0^p(\bar{X}_0)$ with support in
$\supp(\psi_3)$). With these notations,
(\ref{commutdiff0}) is summarized by
\[D^px^\alpha=x^\alpha D^p\]
and we also get
\[D_0^pI_*=I_*D^p, \quad D^pI^*=I^*D_0^p.\]
Using that $\chi_0,\chi_1$ and $[\Delta_{g_0},\psi_2]$
have compact support, we obtain
\begin{equation}\label{lrla}
x^{-\demi}L_R(\la)=\left(D^1+(\Lambda_0-\Lambda)
x^{-1}\chi_0\right)x^\demi
R(\la_0)\chi_1+x^\demi I^*D_0^2R_0(\la)\psi_1I_*,
\end{equation}
\begin{equation}\label{erla}
x^\demi E_R(\la)=x^\demi R_{0R}(\la)+x^\demi\chi_0R(\la_0)\chi_1.
\end{equation}
Similarly, the second indentity of
(\ref{egalpoids}) and the definition of $E_L(\la),L_L(\la)$
imply that
\begin{equation}\label{rxdemi}
\begin{array}{rcl}
R(z)x^\demi&=&I^*\psi_1R_0(z)(I_*x^\demi D^2x^\demi R(z)x^\demi+
\psi_2I_*x^\demi)\\
& &+\chi_1R(\la_0)x^\demi \left(\chi_0+D^1
x^\demi R(z)x^\demi+(Z-\Lambda_0)x^{-\demi}\chi_0R(z)x^\demi\right).
\end{array}
\end{equation}
Combining this last expression with (\ref{lrla}) gives
\[\frac{K(\la,z)}{\Lambda-Z}=\left(D^1+(\Lambda_0-\Lambda)
x^{-1}\chi_0\right)x^\demi R(\la_0)\chi_1R(z)x^\demi\]
\[+x^\demi I^*D_0^2\psi_4R_0(\la)\psi_1I_*x^\demi
x^{-\demi}\chi_1R(\la_0)x^\demi \left(\chi_0+D^1
x^\demi R(z)x^\demi+(Z-\Lambda_0)x^{-\demi}\chi_0R(z)x^\demi\right)\]
\[+ x^\demi I^*D_0^2\psi_4R_0(\la)\psi_1^2R_0(z)\psi_4I_*x^\demi
(D^2x^\demi R(z)x^\demi+\til{\psi}_2).\]
The first two lines of this expression can be extended to
$\{\Re(\la)>\ndemi-\frac{1}{4}\}\cap\{|\Im(\la)|\geq 1\}$
and $\{\Re(z)\geq \ndemi\}\cap\{|\Im(z)|\geq C\}$ as bounded
operators on $L^2(X)$ by using Lemma \ref{casmodel} and
Theorem \ref{cardvod} (to be in the settings of Lemma \ref{casmodel},
we take $\rho:=i_*(x^\demi)\psi_4$ and we remark that
$x^\demi I^*D_0^2\psi_4=I^*D_0^2\rho$ and $\psi_1I_*x^\demi=\rho\psi_1I_*$),
and their $\mc{L}(L^2)$ norms is bounded by
$Ce^{C(|\la|+|z|)}$.

It remains to consider the last line in the expression of
$(\Lambda-Z)^{-1}K(\la,z)$.
Using (\ref{propdiff0}) and Theorem \ref{cardvod}, observe that
\[D^2x^\demi R(z)x^\demi=D^2(\Delta_g+i)^{-1}(D^1x^\demi R(z)x^\demi+x+(i+Z)
x^\demi R(z)x^\demi)\]
can be extended to $\Re(z)=\ndemi$ and $|\Im(z)|\gg 0$ with
$\mc{L}(L^2)$ norm bounded by $Ce^{C|z|}$.
Note that
\begin{eqnarray*}
\rho R_0(\la)\psi_1^2R_0(z)\rho&=&
\rho R_0(\la)\psi_1(\Delta_{g_0}-\Lambda)R_0(\la)R_0(z)
(\Delta_{g_0}-Z)\psi_1R_0(z)\rho\\
&=&\rho(\psi_1+R_0(\la)[\psi_1,\Delta_{g_0}])R_0(\la)R_0(z)
(\psi_1-[\psi_1,\Delta_{g_0}]R_0(z))\rho.
\end{eqnarray*}
Now it is easy to see that it can be rewritten by
\[(\psi_1+\rho R_0(\la)[\psi_1,\Delta_{g_0}]x^{-\demi})
\frac{\rho R_0(\la)\rho -
\rho R_0(z)\rho}{\Lambda-Z}
(\psi_1-x^{-\demi}[\psi_1,\Delta_{g_0}]R_0(z)\rho).\]
Recall that $x^{-\demi}[\psi_1,\Delta_{g_0}]x^{-\demi}
\in\textrm{Diff}^1_0(\bar{X}_0)$ since $[\psi_1,\Delta_{g_0}]$
has compact support in $X_0$.
Using this expression with Lemma \ref{casmodel}, we obtain that
\[x^\demi I^*D_0^2\psi_4R_0(\la)\psi_1^2R_0(z)\psi_4I_*x^\demi
(D^2x^\demi R(z)x^\demi+\til{\psi}_2)\]
extends to $\{\Re(\la)>\ndemi-\frac{1}{4},|\Im(\la)|\geq 1\}$
and $\{\Re(z)\geq \ndemi,|\Im(z)|\geq C\}$ as continuous operators
on $L^2(X)$ with norm  bounded by $Ce^{C(|\la|+|z|)}$.
Fix $z=\ndemi+is$ with $|s|$ large, then all
these estimates prove that
\[||K(\la,z)||_{\mc{L}(L^2)}\leq\left|\Re(\la)-\ndemi\right|Ce^{C|s|}\]
when $\Im(\la)=s$ and $\Re(\la)>\ndemi-\frac{1}{4}$.
Moreover this term is bounded by $\demi$ if
\begin{equation}\label{condition}
\left|\Re(\la)-\ndemi\right|\leq 
\demi C^{-1}e^{-C|\Im(\la)|}, \quad \Im(\la)=s
\end{equation}
Since $K(\la,z)$ is holomorphic in $\la$ in $\{\Re(\la)>\ndemi-\frac{1}{4}
,|\Im(\la)|>1\}$ and $C$ does not depend on $s$, one
can invert $(1+K(\la,z))$ holomorphically if $\la$ satisfies
(\ref{condition}).

Finally the term $K_1(\la,z)$ can be treated in a similar way, using
(\ref{erla}) and (\ref{rxdemi}). So the proof is complete in the general case.
\qed\\

When $g$ is non-trapping we could apply the same method
and prove that there exists a free of resonance region
polynomially close to the critical line. We prefer to only write down 
the  case of constant curvature near infinity since we obtain a 
slightly better result.\\

\textsl{Proof of Theorem \ref{hyper}}: the proof is essentially similar.
When the metric $g$ has constant curvature near
infinity there exist (see \cite{GZ2}) some charts
$(V_j)_{j=1,\dots,M}$ covering a neighbourhood of the boundary
$\pl\bar{X}$ such that each $V_j$ is isometric (note $i_j$ this isometry)
to the open set
\begin{equation}\label{boulehyp}
B:=\{m=(y_1,\dots,y_{n+1})\in\rr^{n+1};y_{n+1}>0, \sum_{i}y_i^2< 1\}
\end{equation}
equipped with the hyperbolic metric
$y^{-2}_{n+1}(\sum_{i}dy_i^2)$. Actually, a parametrix sufficient
for our problem is given in \cite[Prop. 3.1]{GZ2}. Fix
$\la_0\in\{\Re(\la)>\ndemi\}$ and recall that $\Delta_{g_h}$ is the
Laplacian on $\hh^{n+1}$ and $R_h(\la)$ its resolvent studied in
Lemma \ref{bound:uniform}.
Let $\chi_i^j$ be some functions with
support in (\ref{boulehyp}), such that $\chi_1^j=1$ on the support
of $\chi_2^j$ and 
\[\chi_i^j(y_1,\dots,y_{n+1})=\phi_i^j(y_1,\dots,y_n)\psi_i^j(y_{n+1})\]
with $\phi_i^j\in C_0^\infty(\rr^n)$ satisfying
$\sum_ji_j^*\phi^j_2=1$ on $\pl\bar{X}$ and $\psi_i^j\in
C_0^\infty([0,1))$ which is equal to $1$ for $y_{n+1}\leq \delta$
with $\delta<1$ small (see \cite{GZ2} or \cite[Lem. 3.2]{PP} for details).
This implies that 
$\chi:=1-\sum_{j=1}^Mi_j^*\chi^j_2$ has compact support in $X$.
Let $\chi_0\in C_0^\infty(X)$ such that $\chi_0=1$ on the support
of $\chi$. Finally $I_j^*$ and ${I_j}_*$ are defined as in
the proof of previous theorem with our new isometries $V_j\to B$ and the
inclusions $V_j\subset X$, $B\subset\hh^{n+1}$. We then
apply the proof of previous theorem but we replace $E_R,E_L$ by our
new parametrix
\[E_R(\la):=\chi_0R(\la_0)\chi+\sum_{j=1}^MI_j^*\chi_1^jR_h(\la)\chi_2^j{I_j}_*,
\quad E_L(\la):=\chi
R(\la_0)\chi_0+\sum_{j=1}^MI_j^*\chi_2^jR_h(\la)\chi_1^j{I_j}_*\]
and the error terms are
\[L_R(\la):=[\Delta_{g},\chi_0]R(\la_0)\chi+(\Lambda_0-\Lambda)
\chi_0R(\la_0)\chi+\sum_{j=1}^MI_j^*[\Delta_{g_h},\chi_1^j]R_h(\la)\chi_2^j{I_j}_*,\]
\[L_L(\la):=\chi R(\la_0)[\chi_0,\Delta_g]+(\Lambda_0-\Lambda)
\chi R(\la_0)\chi_0+\sum_{j=1}^MI_j^*\chi_2^j
R_h(\la)[\chi_1^j,\Delta_{g_h}]{I_j}_*\]
where $\Lambda:=\la(n-\la)$, $\Lambda_0:=\la_0(n-\la_0)$. As before,
$x_j:={i_j}_*x$ is a boundary defining function of $\{y_{n+1}=0\}$ in $\bar{B}$.
Moreover it is easy to check (see \cite[Prop. 3.1]{GZ2}) that
\begin{equation}\label{propcom}
x_j^{-\demi}[\Delta_{g_h},\chi_i^j]x_j^{-\demi}\in
\textrm{Diff}_0^1(\bbar{\hh}^{n+1}), \quad i=1,2
\end{equation}
where we consider $B\subset \hh^{n+1}$.
We then have the same as (\ref{paramfinal}) with
\begin{eqnarray*}
\frac{K(\la,z)}{\Lambda-Z}&=&x^{-\demi}\left([\Delta_{g},\chi_0]
R(\la_0)\chi+(\Lambda_0-\Lambda)\chi_0 R(\la_0)\chi+
\sum_{j=1}^MI_j^*[\Delta_{g_h},\chi_1^j]R_h(\la)\chi_2^j{I_j}_*\right)\\
& &\x \Big(\chi R(\la_0)\chi_0+\sum_{j=1}^MI_j^*\chi_2^jR_h(z)
(\chi_1^j{I_j}_*+[\Delta_{g_h},\chi_1^j]{I_j}_*R(z))\\
& &\quad\quad\quad\quad\quad\quad+\chi R(\la_0)[\Delta_{g},\chi_0]R(z)-(\Lambda_0-Z)\chi
R(\la_0)\chi_0R(z)\Big)x^{\demi}
\end{eqnarray*}
for $\la,z$ in the physical sheet. Note that Theorem
\ref{cardvod}, Lemma \ref{bound:uniform} and (\ref{propcom}) prove
that all these products extend to $\Re(z)=\ndemi$ and
$\Re(\la)>\ndemi-\frac{1}{8}$  except maybe
\begin{equation}\label{badterm}
x^{-\demi}\left(\sum_{j=1}^MI_j^*[\Delta_{g_h},\chi_1^j]R_h(\la)\chi_2^j{I_j}_*\right)
\left(\sum_{j=1}^MI_j^*\chi_2^jR_h(\la)(\chi_1^j{I_j}_*+
[\Delta_{g_h},\chi_1^j]{I_j}_*R(z))\right)x^\demi.
\end{equation}
Fix $\la_0=\frac{n}{2}+\frac{1}{8}+is_0$ with $s_0>0$ large and
consider some $(\la,z)$ in
\[\mc{O}:=\left\{(\la,z)\in\cc^2; |\la-\la_0|<\frac{1}{4},  |z-\la_0|<
\frac{1}{4}, \Re(z)\geq\ndemi\right\}.\]
Theorem \ref{cardvod} implies that for $\alpha=-\demi,0$
\begin{equation}\label{boundla0z}
\left|\left|x^{-\demi}[\Delta_{g},\chi_0]R(\la_0)\chi
x^{\alpha}+(\Lambda-\Lambda_0)x^{-\demi}\chi_0 R(\la_0)\chi
x^{\alpha}\right|\right|_{\mc{L}(L^2)}\leq C,
\end{equation}
\begin{equation}\label{boundla0z2}
\left|\left|x^{\alpha}\chi R(\la_0)\chi_0x^\demi
+(Z-\Lambda_0)x^{\alpha}\chi R(\la_0)\chi_0R(z)x^\demi +x^{\alpha}\chi
R(\la_0)[\Delta_{g},\chi_0]R(z)x^\demi\right|\right|_{\mc{L}(L^2)}\leq
Cs_0^{-1}
\end{equation}
and $C$ does not depend on $s_0$. Recall that
$I_j^*, {I_j}_*$ are some isometries
$L^2(V_j)\leftrightarrow L^2(B)$ and $x_j={i_j}_*x$, then
Lemma \ref{bound:uniform} and (\ref{propcom}) lead to
\begin{equation}\label{boundrh}
\left|\left|\sum_{j=1}^MI_j^*x_j^{-\demi}[\Delta_{g_h},\chi_1^j]
R_h(\la)\chi_2^jx_j^\demi {I_j}_*\right|\right|
_{\mc{L}(L^2)}\leq C
\end{equation}
\begin{equation}\label{boundrh2}
\left|\left|\sum_{j=1}^MI_j^*x_j^\demi\chi_2^jR_h(z)x_j^\demi
(\chi_1^j{I_j}_*+x_j^{-\demi}[\Delta_{g_h},\chi_1^j]x_j^{-\demi}{I_j}_*
x^\demi R(z)x^\demi)\right|\right|_{\mc{L}(L^2)}
\leq Cs_0^{-1}
\end{equation}
and it remains to consider (\ref{badterm}), which a priori does not exist in
$\mc{O}$. When $V_j\cap V_i=\emptyset$
it is clear that $\chi_2^j{I_j}_*I_i^*\chi_2^i=0$,
so suppose that $V_j\cap V_i\not=\emptyset$. Using that the isometry 
$i_i\circ i_j^{-1}:i_j(V_i\cap V_j)\to i_i(V_i\cap V_j)$ extends to an isometry $I_{ij}$
on $\hh^{n+1}$ and following the proof of
previous theorem, we can see that
\[x_j^\demi R_h(\la)\chi_2^j{I_j}_*I_i^*\chi_2^iR_h(z)x_i^\demi\]
can be expressed for $\la,z$ in the physical sheet $\{\Re(\la)>\ndemi\}$ by
\[(\chi_2^j+x_j^\demi R_h(\la)[\chi_2^j,\Delta_{g_h}]x_j^{-\demi})
\frac{x_j^\demi R_h(\la)\til{x}_j^\demi -x_j^\demi R_h(z)\til{x}_j^\demi}
{\Lambda-Z}(\til{x}_j^{-\demi}[\til{\chi}_2^i,\Delta_{g_h}]R_h(z)\til{x}_j^\demi+\til{\chi}_2^i)I_{ij}^*\]
where $\til{\chi}_2^i:=I_{ij}^*\chi_2^i$, $\til{x}_j:=I_{ij}^*x_i$.
Observe that this operator on $\hh^{n+1}$ extends to $(\la,z)\in\mc{O}$
and satisfies in $\mc{O}$
\begin{equation}\label{boundrlarz}
||x_j^\demi R_h(\la)\chi_2^j{I_j}_*I_i^*\chi_2^iR_h(z)x_i^\demi||_{\mc{L}(L^2,H^1)}
\leq Cs_0^{-1}
\end{equation}
in view of Lemma \ref{bound:uniform} and the following bound (implicitly used
in the proof of previous theorem)
\[\left|\left|\frac{x_j^\demi R_h(\la)\til{x}_j^\demi-x_j^\demi R_h(z)\til{x}_j^\demi}
{\Lambda-Z}\right|\right|_{\mc{L}(L^2,H^1)}\leq
Cs_0^{-1}\sup_{|\la-\la_0|<\frac{1}{4}}||\pl_\la x_j^\demi
R_h(\la)\til{x}_j^\demi||_{\mc{L}(L^2,H^1)}.\]
Combining (\ref{boundla0z})-(\ref{boundrlarz}) with Theorem \ref{cardvod},
one concludes that
\[||K(\la,z)||\leq C|\la-z|\]
for $(\la,z)\in \mc{O}$, hence by taking $z=\ndemi+is_0$ we see
that $1+K(\la,z)$ is holomorphically invertible for
\[\la\in\{\la\in\cc;|\la-z|<\max(C^{-1},\frac{1}{8})\}.\]
Since $C$ does not depend on $s_0$ (thus on $z$), there exists a strip $\{|\Re(\la)-\ndemi|<\eps\}$
where $K(\la,z)$ is invertible except maybe at a finite number of points. 
Finally the term $K_1(\la,z)$ defined as in (\ref{paramfinal})
with our new parametrix can be studied similarly.
\qed

\end{document}